\documentclass[a4paper, BCOR=0.0mm, DIV=calc]{scrartcl}
\usepackage[T1]{fontenc}
\usepackage[utf8]{inputenc}
\usepackage[ngerman]{babel}
\usepackage[osf,sc]{mathpazo}
\usepackage{textcomp}
\usepackage{microtype}
\usepackage{graphicx, color}
\usepackage{amsmath, amssymb, amsfonts}
\usepackage{booktabs}
\usepackage{nicefrac}
\usepackage{wasysym,pifont}
\usepackage{tikz}
\usepackage{graphicx}
\usepackage{float}
\usepackage{wrapfig}
\usetikzlibrary{calc,intersections}
\KOMAoptions{twoside=false, twocolumn=false, headinclude=false, footinclude=false, mpinclude=false, pagesize=auto}
\recalctypearea

\setkomafont{section}{\mdseries\scshape\Large}
\addtokomafont{subsection}{\mdseries\scshape}
\setkomafont{title}{\mdseries\scshape}

\begin{document}
\title{On Solids whose surface can be unfolded onto a plane
\footnote{
Original title: "' De solidis quorum superficiem in planum explicare licet"', first published in "`\textit{Novi Commentarii academiae scientiarum Petropolitanae} (1771), 1772, pp. 3-34"', reprinted in  "`\textit{Opera Omnia}: Series 1, Volume 28, pp. 161 - 186 "', Eneström-Number E419, translated by: Alexander Aycock, for the "`Euler-Kreis Mainz"'.}}
\author{Leonhard Euler}
\date{ }
\maketitle

\paragraph*{§1} The property of the cylinder and the cone that their surface can be unfolded onto a plane is very well-known, and this property is even extended to all cylindrical and conic bodies whose bases have any arbitrary shape; the sphere on the other hand does not enjoy this property, since its surface  can not be unfolded onto a plane  by any means and it can not be covered by a planar surface; this gives rise to the curious and intriguing question, whether except for cones and cylinders other classes of solids exist whose surface can be unfolded onto the plane in the same way or not. Therefore, I decided to consider the following problem in this dissertation:\\[2mm]
\textit{To find a general equation for all solids whose surface can be unfolded onto a plane},\\[2mm]
whose solution I will attempt to give in various ways. \newpage

\subsection*{First Solution derived from mere analytical Principles}

\paragraph*{§2} 

Let (Fig. 1) $Z$ be an arbitrary point on the surface of the solid in question; let the location of this point, as it customary now, be expressed by three mutually orthogonal coordinates $AX=x$, $XY=y$ and $YZ=z$ such that an equation between these three coordinates is to be found, by means of which the problem will then be solved.
\begin{center}
\includegraphics[scale=0.4]{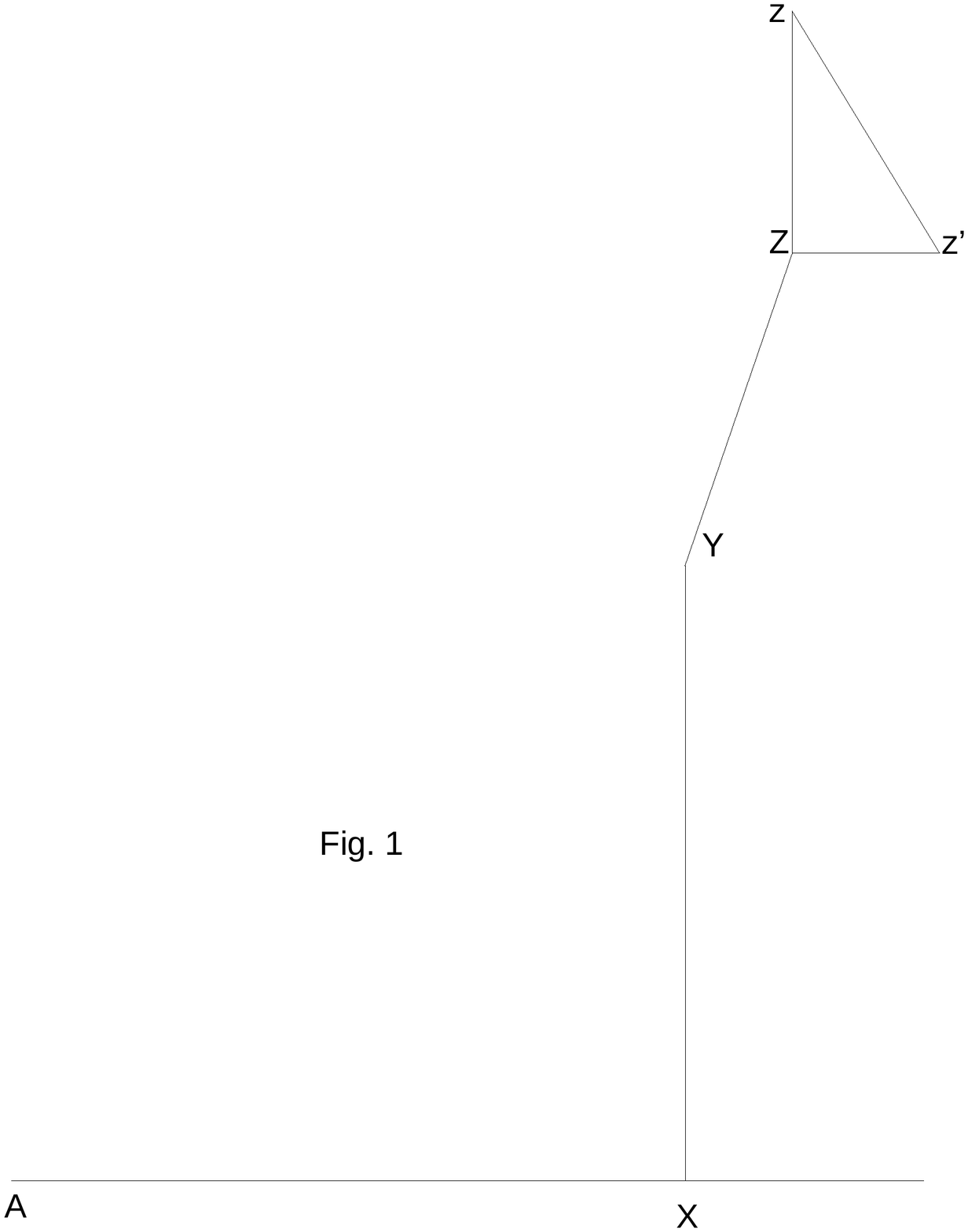}
\end{center}
Further, let us assume that the surface of a solid of such a kind has already been unfolded onto the plane and it is represented in figure 2, in which the point $Z$ falls on $V$; next define  location of this point $V$  by two orthogonal coordinates in such a way that $OT=t$ and $TV=u$. Then it is manifest that the first three coordinates $x$, $y$ and $z$ have to depend on these two $t$ and $u$ in a certain way, and hence every single one of them can be considered a certain function of $t$ and $u$.

\paragraph*{§3} 

In order to introduce this condition into the calculation in a more convenient way let us consider it in terms of differentials and, since $x$, $y$ and $z$ are functions of the two variables $t$ and $u$, let us define their differentials by  these formulas:

\begin{equation*}
dx= ldt+\lambda du, \quad dy=mdt+\mu du \quad \text{and} \quad dz= ndt+\nu du;
\end{equation*}
here, since the letters $l$, $m$, $n$ and $\lambda$, $\mu$, $\nu$ in the same way denote certain functions of the two variables $t$ and $u$, it is clear from the nature of functions of this kind\footnote{Euler means functions of two and more variables and refers to what we nowadays know as Schwarz's theorem for the partial derivatives of functions of two or more variables.} that it has to be:

\begin{equation*}
\left(\frac{dl}{du}\right)=\left(\frac{d \lambda}{dt}\right),\quad \left(\frac{dm}{du}\right)=\left(\frac{d \mu}{dt}\right) \quad \text{and} \quad \left(\frac{dn}{du}\right)=\left(\frac{d \nu}{dt}\right).
\end{equation*}

\paragraph*{§4}
Now, in the unfolded surface (Fig.2) let us,
\begin{center}
\includegraphics[scale=0.35]{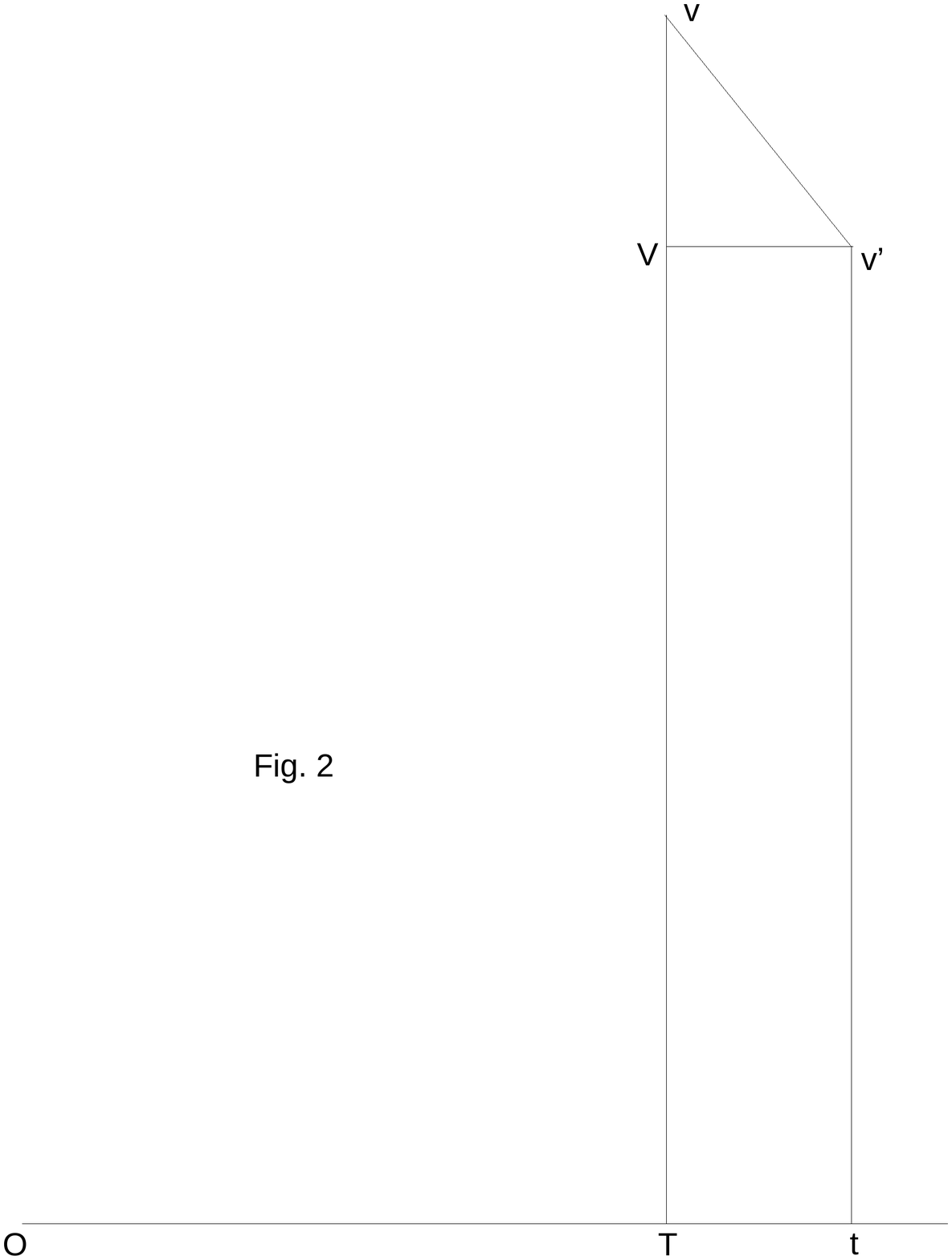}
\end{center} 
 except for the point $V$, contemplate two other points infinitely close to it\footnote{Euler refers to infinitesimal quantities several times in this paper; for further explanation of this concept see his book “Institutiones calculi differentialis” (E212).}, $v$ and $v^{\prime}$; for the first of them let the coordinates be
 
 \begin{equation*}
 OT= t \quad \text{and} \quad Tv=u+du,
\end{equation*}  
for the latter on the other hand:

\begin{equation*}
Ot=t+dt \quad \text{and} \quad tv^{\prime}=u,
\end{equation*}
such that the points $V$ and $v$ have the common abscissa $OT=t$, but the points $V$ and $v^{\prime}$ have the common ordinate $=u$. Having drawn the infinitely short lines $Vv^{\prime}$ and $vv^{\prime}$ the sides of the elementary triangle $Vvv^{\prime}$ are determined in such a way that it is:

\begin{equation*}
Vv=du,\quad Vv^{\prime}=dt \quad \text{and} \quad vv^{\prime}=\sqrt{du^2+dt^2},
\end{equation*}
and now it is easily understood that the same triangle has to be found also in the surface of the solid in question.

\paragraph*{§5}

Therefore, in the surface of the solid let $z$ and $z^{\prime}$ be the points corresponding to the points $v$ and $v^{\prime}$  and let us see how the three coordinates behave for those points $z$ and $z^{\prime}$. But, the way how the point $Z$ is defined by these three coordinates, the first $=x$, the second $=y$, the third $=z$, which are all functions of the two variables $t$ and $u$, since for the point $v$ the abscissa $t$ remains the same, but the ordinate $u$ on the other hand is augmented by its differential $du$, the three coordinates for the point $z$ of the solid will behave as this:

\begin{equation*}
\text{I.} \quad x+\lambda du, \quad \text{II.} \quad y+\mu du \quad \text{and} \quad \text{III.} \quad z+\nu du; 
\end{equation*}
in like manner, because for the point $v^{\prime}$ the ordinate $u$ remains the same, the abscissa $t$ on the other hand is augmented by its differential $dt$, the three coordinates for the point $z^{\prime}$ will be:

\begin{equation*}
\text{I.}\quad x+ldt, \quad \text{II.} \quad y+mdt \quad \text{and} \quad \text{III.} \quad z+ndt.
\end{equation*}

\paragraph*{§6}

But it is known, if for any arbitrary point on the surface of a solid the coordinates were $x$, $y$ and $z$, but for another infinitely close point the coordinates were $x^{\prime}$, $y^{\prime}$ and $z^{\prime}$, that then the distance of the points will be:

\begin{equation*}
=\sqrt{(x-x^{\prime})^2+(y-y^{\prime})^2+(z-z^{\prime})^2};
\end{equation*}
hence, we will have for the single sides of the triangle $Zzz^{\prime}$:

\begin{alignat*}{9}
&1^{\circ} \quad && Zz&&= du \sqrt{\lambda^2 +\mu^2 +\nu^2},\\
&2^{\circ} \quad && Zz^{\prime}&&= dt \sqrt{l^2 +m^2 +n^2}
\end{alignat*}
and
\begin{equation*}
3^{\circ} \quad zz^{\prime}=\sqrt{(\lambda du-l dt)^2+(\mu du -mdt)^2+(\nu du -ndt)^2}
\end{equation*}
or
\begin{equation*}
zz^{\prime}=\sqrt{dt^2(ll+mm+nn)+du^2(\lambda \lambda +\mu \mu +\nu \nu)-2dtdu(l \lambda +m \mu+n \nu)}.
\end{equation*}

\paragraph*{§7}

Now, because the surface of the solid has to  agree completely with the planar figure (Fig. 2), it is necessary that the triangles $Zzz^{\prime}$ and $Vvv^{\prime}$ are not only equal but also similar and hence the sides  equal homologues\footnote{By this Euler means that that triangles must actually be congruent; confer with the following equations.}, namely:

\begin{equation*}
\text{I}^{\circ}. \quad Zz=Vv, \quad \text{II}^{\circ}. \quad Zz^{\prime}=Vv^{\prime} \quad \text{and} \quad \text{III}^{\circ}. \quad zz^{\prime}=vv^{\prime}, 
\end{equation*}
whence we obtain the following equations:
\begin{alignat*}{9}
&\text{I}^{\circ}. \quad && \lambda^2 +\mu^2 +\nu^2 =1 \\
&\text{II}^{\circ}. \quad && l^2 +m^2 +n^2 =1 \\
&\text{III}^{\circ}. \quad && dt^2(l^2+m^2+n^2)+du^2(\lambda^2 +\mu^2 +\nu^2)-2 dt du (l \lambda +m \mu +n \nu)=dt^2+du^2;
\end{alignat*}
the  third equation, using the first two, is reduced to this one:

\begin{equation*}
l \lambda + m \mu +n \nu =0;
\end{equation*}
these  three equations contain the solution of our problem, whence it is understood that it can be reduced to the following analytical problem:\\[2mm]

\subsection*{Analytical Problem}

\textit{Given two variables $t$ and $u$, to find six functions $l$, $m$, $n$ and $\lambda$, $\mu$, $\nu$ of such a kind that the following six conditions are fulfilled:}

\begin{alignat*}{9}
&\text{I}^{\circ}. \quad \left(\frac{dl}{du} \right)=\left(\frac{d \lambda}{dt}\right), \quad \text{II}^{\circ}. \quad \left(\frac{dm}{du}\right)=\left(\frac{d \mu}{dt}\right), \quad \text{III}^{\circ}, \quad \left(\frac{dn}{du}\right)=\left(\frac{d \nu}{dt}\right), \\
&\text{IV}^{\circ}. \quad ll+mm+nn=1,\quad \text{V}^{\circ}. \quad \lambda \lambda +\mu \mu +\nu \nu =1, \\
&\text{VI}^{\circ} \quad l\lambda +m \mu +n \nu =0; 
\end{alignat*}
this problem, considered in itself, seems to be very difficult;  we will nevertheless be able to exhibit a beautiful solution of it below. \newpage

\subsection*{Second Solution derived from geometrical Principles}

\paragraph*{§8}

In order to derive this solution from first principles let us consider either prismatic or pyramidal bodies which, having excluded the bases, are understood to be covered by a chart, and on this chart rectilinear sharp bends will then be detected which are either parallel to each other or converge to a certain point, the vertex of the pyramid, of course. Denote these, either parallel or converging, straight lines,  by the letters $Aa$, $Bb$, $Cc$, $Dd$ etc. Therefore, if the chart is unfolded onto the plane, in it the same straight lines  $Aa$, $Bb$, $Cc$ etc. will occur, and they will be either parallel to each other or converge to a certain point. Hence vice versa, if  such straight lines are drawn on the planar chart, according to which the chart can be folded, it will be apt for covering a certain prismatic or pyramidal body. 

\paragraph*{§9}

It will even be possible to draw the lines $Aa$, $Bb$, $Cc$, $Dd$ etc. on the chart arbitrarily such that they are neither parallel to each other nor converge to a certain point, as long as they never cross each other, as figure 3 shows; 
\begin{center}
\includegraphics[scale=0.3]{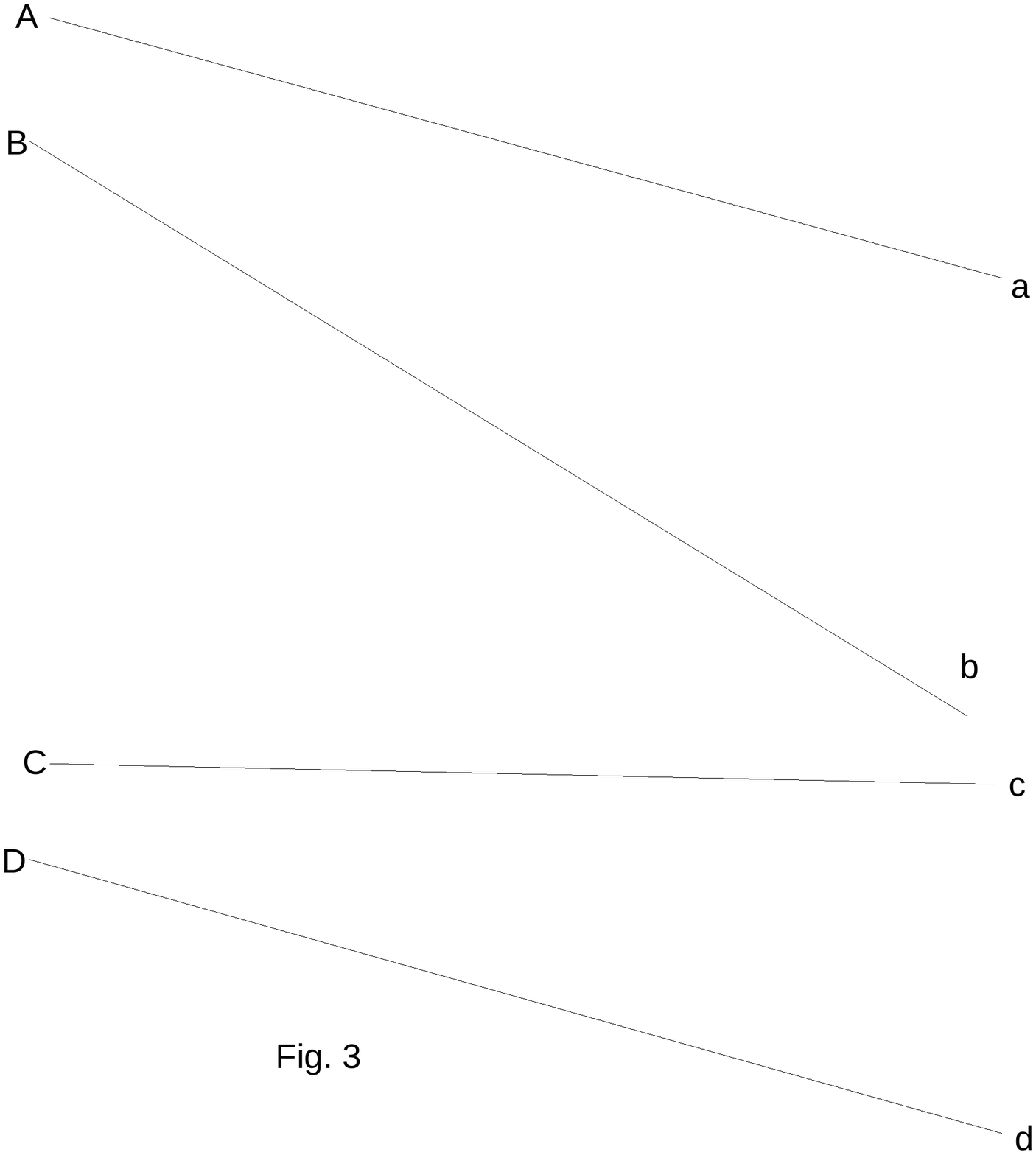}
\end{center} 
for, no matter how this chart is folded according to these lines, it will always possible to conceive a solid of such a kind, which would then by covered by this folded chart. From this it is  plain that beyond prismatic or pyramidal bodies other classes of bodies exist which can be covered by a chart this way and whose surface can therefore be unfolded onto a plane.

\paragraph*{§10}

Therefore, on the surface of these bodies any number of straight lines $Aa$, $Bb$, $Cc$, $Dd$ etc. will be given which, even though they are neither parallel nor converge to a certain point, will nevertheless be of such a nature that any two, very close\footnote{By this Euler means infinitely close.} to each other, such as $Aa$ and $Bb$ or $Bb$ and $Cc$ or $Cc$ and $Dd$ etc., if they are not parallel, at least intersect in a single point if they were elongated; for, if this would not happen, the space between the two lines intersected with the body would not be planar and therefore  it would not be possible to unfold the surface onto the plane, although  there are arbitrarily many straight lines $Aa$, $Bb$, $Cc$ etc. in it. From this we conclude that it does not suffice for bodies, meeting the condition that they can be unfolded onto the plane,  that it is possible to draw an arbitrary number of lines $Aa$, $Bb$, $Cc$ etc. on them; but furthermore it is required that two lines infinitely close to each other are in the same plane and the space contained between them is itself planar.

\paragraph*{§11}

Now, let us increase the number of the straight lines $Aa$, $Bb$, $Cc$ etc. to infinity such that our body obtains a curved surface everywhere, as our problem postulates it because of the law of continuity.  And now it is indeed immediately clear that a surface of such a kind has to be of such a nature that from any arbitrary point assumed in it at least one straight line can be drawn which lies in the surface completely; but this condition alone does not exhaust the whole character of our problem, but additionally it is necessary that any two lines of this kind, infinitely close to each other, lie in the same plane; this means that, if they are not parallel, they at least meet in one point, if one would elongate them. Hence, if those single lines are elongated to the point of intersection  this way, all these points of intersection will be found to lie on a certain curve. Because this curve does not lie in one plane completely, it will have two curvatures and will be of such a nature that its single elements, if they are elongated, exhibit those lines $Aa$, $Bb$, $Cc$ etc. mentioned above  in the surface of the body.

\paragraph*{§12}

Therefore, as any body convenient for our problem leads to a certain curve with two curvatures, so vice versa having assumed a curve of this kind arbitrarily we will be able to determine a body from it, which solves our problem. But first project such a curve onto the plane of the diagram, and let (Fig. 4) its projection be $aUu$;

\begin{center}
\includegraphics[scale=0.5]{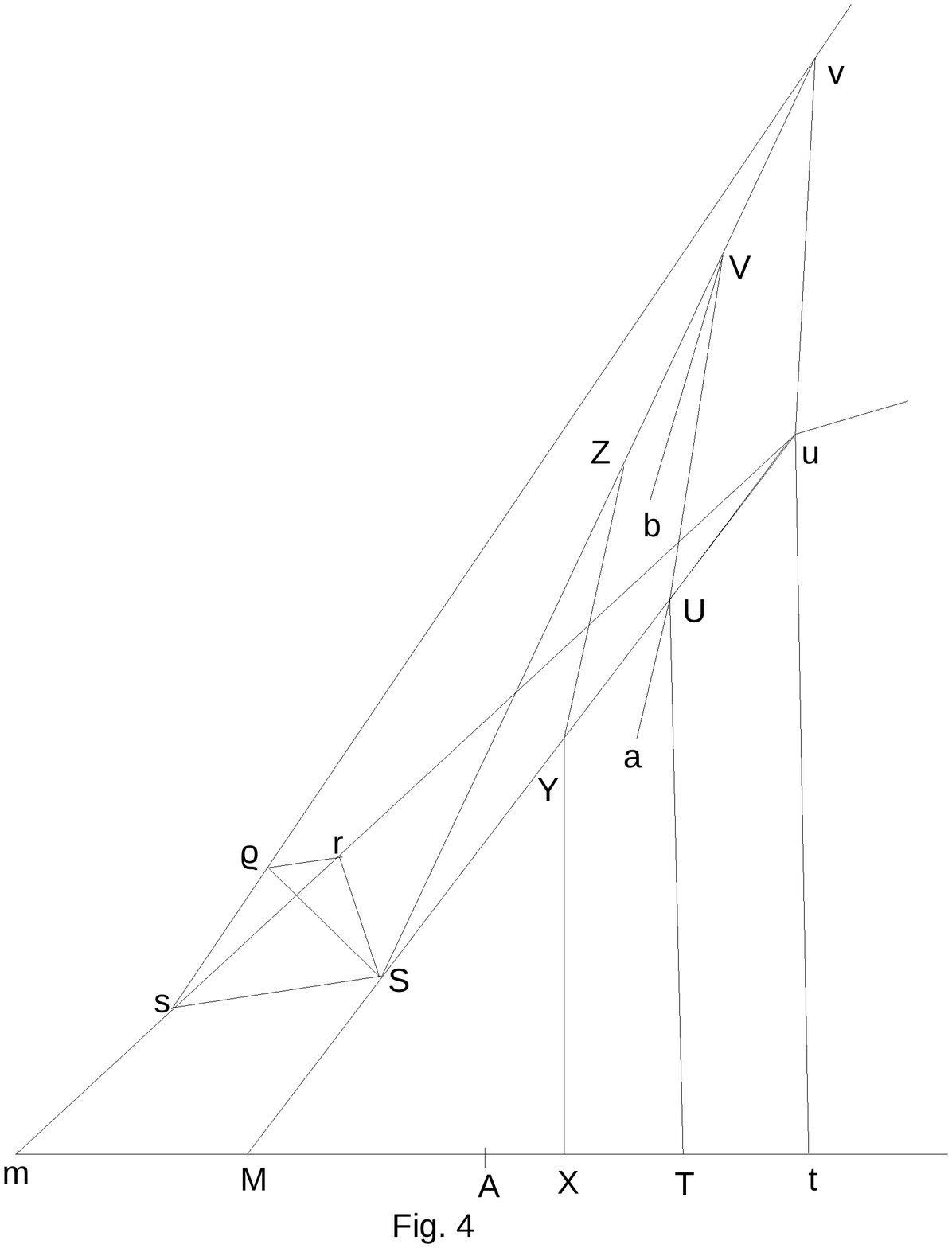}
\end{center} 
for this we want to put the abscissa $AT=t$ and the ordinate $Tu=u$ such that an equation between $t$ and $u$ is considered as given, and let $UM$ be the tangent of this curve at the point $U$, the line $um$ on the other hand the tangent at point $u$ infinitely close to $U$; having constituted all this let $bVv$ be the curve with the two curvatures, whose ordinate orthogonal to our plane shall be put $UV=v$, and let $v$ be the closest  point  on the same curve, and  draw tangents starting from both points $V$, $v$; let  the first of these tangents, $VS$,  intersect the line $UM$ in the point $S$ and let the other, $vs$,  intersect the line $rm$ in the point $s$. Here, we certainly could have drawn the infinitely close tangents in the points $u$ and $v$, but, since it will be necessary in the following, it seemed advisable to indicate them in the figure already, as we did here.

\paragraph*{§13}

Therefore, because the nature of the curve $bVv$ is expressed by two equations between the coordinates $AT=t$, $Tu=u$ and $UV=v$, the letters $u$ and $v$ can be considered as a function of $t$, whence at the same time the position of both tangents $UM$ and $VS$ will be defined; hence let us  to call the angles $TUM= \zeta$ and $UVS= \vartheta$; and having put the element $Tt=dt$ it will be

\begin{equation*}
du= \frac{dt}{\tan \zeta}, \quad Uu= \frac{dt}{\sin \zeta},
\end{equation*}
then on the other hand it will be

\begin{equation*}
dv= \frac{dt}{\sin \zeta \tan \vartheta}
\end{equation*}
and finally the element of the curve will be

\begin{equation*}
Vv= \frac{dt}{\sin \zeta \sin \vartheta}.
\end{equation*}
But, for the position of the tangents we will have

\begin{equation*}
TM= u \tan \zeta, \quad UM= \frac{u}{\cos \zeta},
\end{equation*}
the line on the other hand

\begin{equation*}
US= v \tan \vartheta \quad \text{and} \quad VS= v \sec \vartheta =\frac{v}{\cos \vartheta}.
\end{equation*}

\paragraph*{§14}

Since now the whole line $VS$ lies in the surface of the body in question, let us take an indefinite point $Z$ on that curve, whence having dropped the perpendicular $ZY$ to the plane of the diagram from that point $Z$ and having drawn the normal $YX$ from the point $Y$ to the axis $AT$,  we will have the three coordinates we contemplated above, of course $AX=x$, $XY=y$ and $YZ=z$, for the surface in question; therefore, the correct equation between them, by which the nature of this surface is expressed is to be investigated.

\paragraph*{§15}

For this aim, let us call the indefinite interval $VZ=s$ which therefore is a variable quantity not depending on the point $V$ and hence is to be  distinguished carefully from the variable $t$; recall that not only the two ordinates $TU=u$ and $UV=v$ are functions of this variable $t$, but also the two angles $\zeta$ and $\vartheta$ are. Hence, we obtain

\begin{equation*}
ZY=z=v-s \cos \vartheta
\end{equation*}
and the interval

\begin{equation*}
UY= s \sin \vartheta,
\end{equation*}
whence we further conclude

\begin{equation*}
XY= y= u-s \sin \vartheta \cos \zeta
\end{equation*}
and

\begin{equation*}
XT= s \sin \vartheta \sin \zeta,
\end{equation*}
and so we finally obtain the abscissa

\begin{equation*}
AX= x= t-s \sin \vartheta \sin \zeta,
\end{equation*}
such that by means of the two variables $t$ and $s$ our three coordinates are determined succinctly this way:

\begin{alignat*}{9}
&\text{I}^{\circ}. \quad && x &&= t- s \sin \vartheta \sin \zeta, \\
&\text{II}^{\circ}. \quad && y &&=u-s \sin \vartheta \cos \zeta, \\
&\text{III}^{\circ}. \quad && z &&=v- s \cos \vartheta.
\end{alignat*}

\paragraph*{§16}

Therefore, against all expectations it happens here that we even found algebraic formulas for the three coordinates $x$, $y$, $z$, if one takes algebraic functions of $t$ for the quantities $u$ and $v$. For, these functions are completely arbitrary; but having assumed them, the two angles $\zeta$ and $\vartheta$ are determined in such a way that $\tan \zeta = \frac{dt}{du}$ or

\begin{equation*}
\sin \zeta = \frac{dt}{\sqrt{dt^2+du^2}} \quad \text{and} \quad \cos \zeta =\frac{du}{\sqrt{dt^2+du^2}},
\end{equation*}
then on the other hand

\begin{equation*}
\tan \vartheta =\frac{dt}{dv \sin \zeta}=\frac{\sqrt{dt^2+du^2}}{dv}
\end{equation*}
and hence

\begin{equation*}
\sin \vartheta =\frac{\sqrt{dt^2+du^2}}{\sqrt{dt^2+du^2+dv^2}} \quad \text{and} \quad \cos \vartheta = \frac{dv}{\sqrt{dt^2+du^2+dv^2}}.
\end{equation*}
But therefore, if vice versa the two angles $\zeta$ and $\vartheta$ were given in terms of the variables $t$, the ordinates $u$ and $v$ will be found expressed by the following integral formulas

\begin{equation*}
u= \int \frac{dt}{\tan \zeta} \quad \text{and} \quad v= \int \frac{dt}{\sin \zeta \tan \vartheta}.
\end{equation*}

\paragraph*{§17}

Therefore,  all solids whose surface can be unfolded onto the plane are necessarily contained in these formulas. Therefore, it will especially be worth one's while to show how the conic bodies are contained in them, since the cylindrical bodies are seen  to be already contained in the conic bodies by moving the vertex to infinity. Therefore, let the point $V$ be the vertex of the cone; since this vertex is fixed, the variables $t$, $u$ and $v$ will also have constant values. Since there is no obstruction that this vertex is taken in the fixed point $A$ itself, we will be able to put $t=0$, $u=0$ and $v=0$; but then because of

\begin{equation*}
\tan \zeta = \frac{dt}{du} \quad \text{and} \quad \tan \vartheta = \frac{dt}{dv \sin \zeta}=\frac{\sqrt{dt^2+du^2}}{dv}
\end{equation*}
these angles $\zeta$ and $\vartheta$ turn out to be indefinite; nevertheless they are indefinite only in such a way that the one can be considered as a certain function of the other, since all things extending to the position of the lines $VS$ are to be referred to one single variable.

\paragraph*{§18}

Therefore, because it is $t=0$, $u=0$ and $v=0$, we will have:

\begin{alignat*}{9}
& &&\text{I}^{\circ}. \quad && x &&=- s \sin \vartheta \sin \zeta, \\
& &&\text{II}^{\circ}. \quad && y &&=- s \sin \vartheta \cos \zeta, \\
&\text{and} \quad &&\text{III}^{\circ}. \quad && z &&=- s \cos \vartheta,
\end{alignat*}
whence it is
\begin{equation*}
\frac{x}{y}=\tan \zeta \quad \text{and} \quad \frac{x}{z}= \tan \vartheta \sin \zeta,
\end{equation*}
from which it is concluded

\begin{equation*}
\sin \zeta = \frac{x}{\sqrt{xx+yy}}
\end{equation*}
and hence from this it follows

\begin{equation*}
\tan \vartheta =\frac{\sqrt{xx+yy}}{z};
\end{equation*}
therefore, since $\tan \vartheta$  becomes equal to an arbitrary function of $\tan \zeta$, we will have such an equation:

\begin{equation*}
\frac{\sqrt{xx+yy}}{z}= \Phi : \left(\frac{x}{y}\right),
\end{equation*}
and so the quantity $\frac{\sqrt{xx+yy}}{z}$ will become equal to a homogeneous function of no dimension of $x$ and $y$\footnote{Euler means that this is a homogeneous function with degree of homogeneity $=0$. The term degree of homogeneity did not exist back then, of course.} and hence further the quantity $z$ will become equal to a homogeneous function of one dimension of $x$ and $y$, or, what reduces to the same, the equation between $x$, $y$ and $z$ will be of such a nature that in it the three variables $x$, $y$ and $z$ will add up to the same number of dimensions everywhere. Therefore, if one of the coordinates $x$, $y$ and $z$ becomes infinite, the equation for the solid will only contain the  two remaining variables which is a criterion for cylindrical bodies.

\paragraph*{§19}

We do not spend more time on the consideration of other solids solving our problem here, because below we will explain a third method, by which we are able to conceive of all species of bodies of this kind a lot more easily. While this second method provided us with such a simple solution, although  by means of the first method  hardly any solution could be hoped for, we will nevertheless now also be able to expand the first solution further and even resolve those analytical formulas, what on first sight seemed to be exceedingly difficult; having done this it will illustrate the whole analysis very well. To do this it  will only be necessary that we carefully reduce this second solution to the elements of the first. \newpage

\subsection*{Application of the second Method to the first Solution}

\paragraph*{§20}

Since in the second solution we have already found formulas for the three coordinates $x$, $y$ and $z$, in which  the nature of the solid is contained, we will have to elaborate on this a bit more in order to  investigate  formulas for the planar figure onto which the surface of the solid is unfolded. Here,  especially that curve $bVv$ with the two curvatures is to be studied more accurately, which by  unfolding  the surface is also reduced to the plane. But because this curve by means of inflections  can be reduced to the plane in infinitely many ways and can even be stretched out into a straight line this way, it is especially to be inquired, how this reduction to the plane has to be actually done. From the things mentioned above it is indeed manifest that this reduction has to happen in such a way that (Fig. 4) any two infinitely close tangents $VS$ and $vs$ conserve the same mutual position to each other or that the angle $Svs$ enclosed between them remains the same. Of course, the curve $bVv$ itself is to be reduced to the plane in such a way that any two infinitely close  elements of it conserve the same inclination to each other.

\paragraph*{§21}

Therefore, the the main task reduces to finding the infinitely small angle $Svs$; for this aim one has to start from the angle $MUm$. But because it is

\begin{equation*}
\text{angle } TUM=\zeta \quad \text{and} \quad \text{angle }tum=\zeta +d \zeta,
\end{equation*}
it manifestly follows that the angle $Mum=d\zeta$; further, because we already found $US= v \tan \vartheta$ above,  from the nature of differentials it will be:

\begin{equation*}
us= v \tan  \vartheta +d(v \tan \vartheta)= v \tan \vartheta +dv \tan \vartheta +\frac{v d \vartheta}{\cos^2 \vartheta},
\end{equation*}
where

\begin{equation*}
dv= \frac{dt}{\sin \zeta \tan \vartheta};
\end{equation*}
 therefore, because it is

\begin{equation*}
Uu= \frac{dt}{\sin \zeta},
\end{equation*}
it will be

\begin{equation*}
Us= v \tan \vartheta +dv \tan \vartheta +\frac{v d \vartheta}{\cos^2 \vartheta}-\frac{dt}{\sin \zeta}= v \tan \vartheta +\frac{v d \vartheta}{\cos^2 \vartheta}.
\end{equation*}
Therefore, from $S$ drop the perpendicular $Sr$ to $Us$ so that one has

\begin{equation*}
rs= \frac{v d\vartheta}{\cos^2 \vartheta},
\end{equation*}
then it will indeed be

\begin{equation*}
Sr= v d \zeta \tan \vartheta,
\end{equation*}
whence also the element $Ss$ could be defined, if it would be of any necessity.

\paragraph*{§22}

Now, from the point $r$ let us drop the perpendicular $r\rho$ to the tangent $vs$, that having drawn  $S\rho$ it becomes normal to $vs$, where it is to be noted that the triangle $Sr\rho$ will have an right angle at $r$, because $Sr$ is normal to the plane $sUV$ itself. Because 

\begin{equation*}
\text{angle }rs \rho =90^{\circ}-\vartheta,
\end{equation*}
it will be

\begin{equation*}
r \rho =sr \cdot \sin r s \rho = \frac{v d \vartheta}{\cos \vartheta},
\end{equation*}
whence it is calculated

\begin{equation*}
S \rho = \sqrt{vv d\zeta^2 \tan^2 \vartheta +\frac{vvd\vartheta^2}{\cos^2 \vartheta}}=\frac{v}{\cos \vartheta}\sqrt{d \zeta^2 \sin^2 \zeta +d \vartheta^2}.
\end{equation*}
Therefore, because it is $VS= \frac{v}{\cos \vartheta}$, hence it is concluded

\begin{equation*}
\text{angle }SVs= \frac{S \rho}{VS}= \sqrt{d \zeta^2 \sin^2 \vartheta +d \vartheta^2}.
\end{equation*}

\paragraph*{§23}

Therefore, we found the angle $SVs$ which the two infinitely close elements of the curve enclose; from this one can determine the radius of curvature of this curve in the point $V$ very quickly; of course, it is

\begin{equation*}
\frac{Vv}{SVs}= \frac{dt}{\sin \zeta \sin \vartheta \sqrt{d\zeta^2\sin^2 \vartheta +d\vartheta^2}},
\end{equation*}
which task is therefore not impeded by the two curvatures; it is enough to have remembered this in the transition. But because here the main issue is  the determination of the elementary angle $SVs$, let us call the angle $SVs=d \omega$ such that it is

\begin{equation*}
d \omega =\sqrt{d \zeta^2 \sin^2 \vartheta +d \vartheta^2} \quad \text{or} \quad d\omega^2 -d \vartheta^2=d \zeta^2 \sin^2 \vartheta;
\end{equation*}
here, because the two angles $\zeta$ and $\vartheta$ are determined by the variable $t$ and  the two ordinates $u$ and $v$ are also functions of $t$, it is clear that also the angle $\omega$ has to be considered as a function of the same variable $t$.

\paragraph*{§24}

Now, according to the prescriptions given above (Fig. 5)
\begin{center}
\includegraphics[scale=0.4]{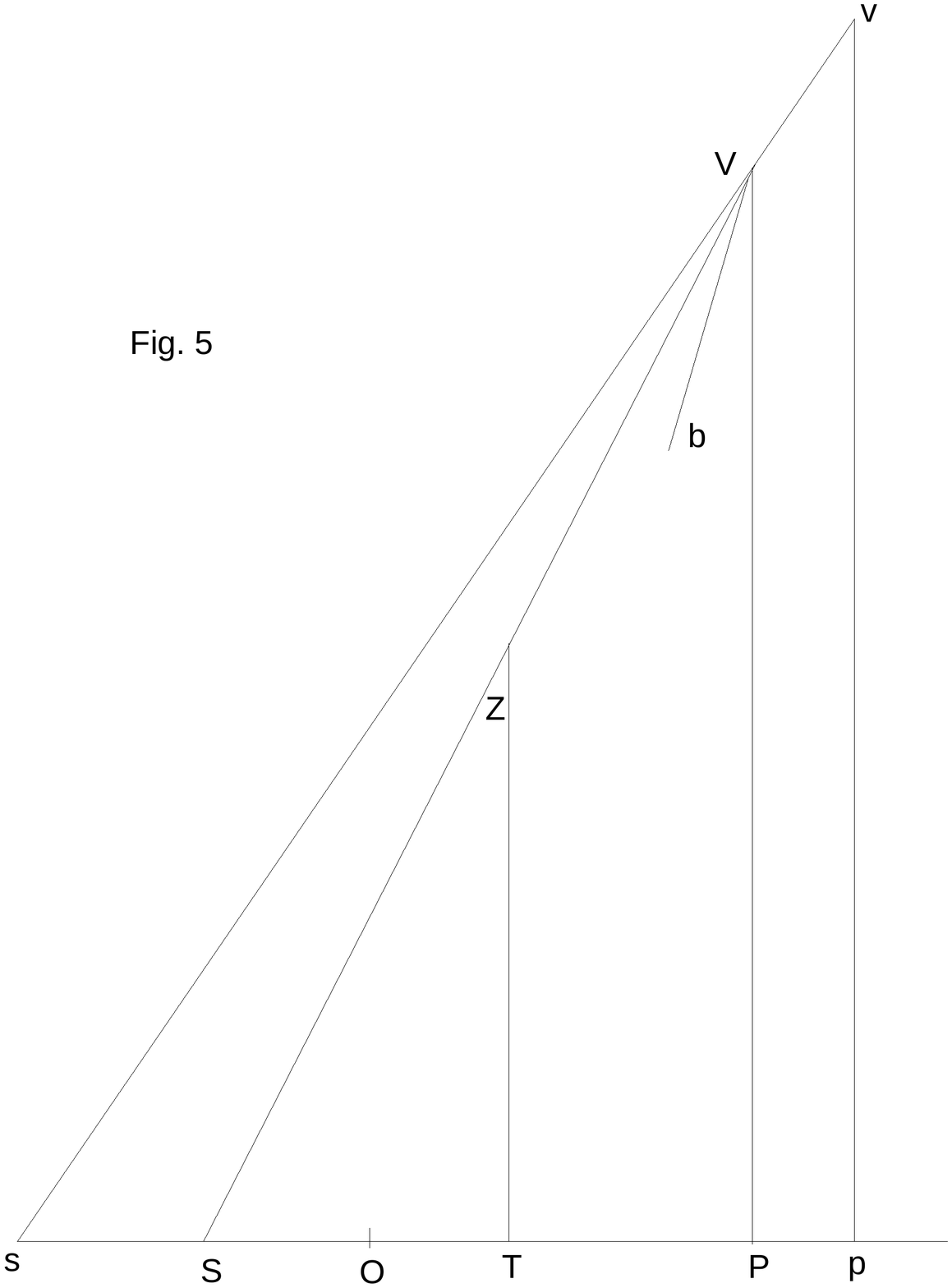}
\end{center}
 let the curve $bVv$ with the two curvatures be described in the plane, such that the angle $SVs$ between two tangents infinitely close to each other will be $=d \omega$, and having related this curve to the axis $OP$ by means of the ordinate $PV$ it is evident that the angle $PVS$ will be $=\omega$. But let us put these coordinates $OP=p$ and $PV=q$, and we will have

\begin{equation*}
\frac{dp}{dq}=\tan \omega,
\end{equation*}
and the element of the curve is

\begin{equation*}
Vv= \frac{dp}{\sin \omega},
\end{equation*}
but on the other hand by means of the preceding coordinates $t$, $u$ and $v$ with the angle $\zeta$ and $\vartheta$ the same element was

\begin{equation*}
Vv= \frac{dt}{\sin \zeta \sin \vartheta},
\end{equation*}
whence  as a logical consequence we obtain

\begin{equation*}
dt \sin \omega = dp \sin \zeta \sin \vartheta,
\end{equation*}
which combined with the equation $\frac{dp}{dq}=\tan \omega$ will give the following integral values for the present coordinates $p$ and $q$

\begin{equation*}
p= \int \frac{dt \sin \omega}{\sin \zeta \sin \vartheta} \quad \text{and} \quad q= \int \frac{dt \cos \omega}{\sin \zeta \sin \vartheta};
\end{equation*}
having found these quantities $p$ and $q$, which likewise are functions of the same variable $t$, take the interval $VZ=s$, which is the other variable to be introduced into the calculation, and having dropped the perpendicular $ZT$ from the point $Z$ to the axis, we find

\begin{equation*}
OT=p-s \sin \omega \quad \text{and} \quad TZ=q- s \cos \omega.
\end{equation*}

\paragraph*{§25}

Therefore, since we obtained the determination for the point $Z$ reduced to the plane, let us put its coordinates $OT=T$ and $TZ=U$, which are defined by  the two variables $t$ and $s$ in such a way  that it is

\begin{alignat*}{9}
&T&&=p&&-s \sin \omega &&= \int \frac{dt \sin \omega}{\sin \zeta \sin \vartheta}&&-s \sin \omega, \\
&U&&=q&&-s \cos\omega &&= \int \frac{dt \cos \omega}{\sin \zeta \sin \vartheta}&&-s \cos \omega,
\end{alignat*}
where it is to be noted that the angle $\omega$ depends on the angles $\zeta$ and $\vartheta$ in such a way that it is

\begin{equation*}
d \omega = \sqrt{d\zeta^2 \sin^2 \vartheta +d\vartheta^2}.
\end{equation*}
These coordinates $T$ and $U$ are indeed the same  we  denoted by the letters $t$ and $u$ in the first solution; hence having made the same change there the formulas found  for the solid there reduce to these 

\begin{equation*}
dx= l dT +\lambda dU, \quad dy=mDT+\mu dU, \quad dz=ndT+\nu dU
\end{equation*}
while the conditions, we found there, remain the same, of course:

\begin{equation*}
ll+mm+nn=1, \quad \lambda \lambda +\mu \mu+\nu \nu=1, \quad \text{and} \quad l \lambda +m \mu +n \nu=0.
\end{equation*}
But here  for the same coordinates $x$, $y$ and $z$ for the solid we found the following values:

\begin{equation*}
x=t-s \sin \vartheta \sin \zeta, \quad y= u-s \sin \vartheta \cos \zeta \quad \text{and} \quad z=v-s \cos \vartheta,
\end{equation*}
which because of

\begin{equation*}
du=\frac{dt}{\tan \zeta} \quad \text{and} \quad dv=\frac{dt}{\sin \zeta \tan \vartheta},
\end{equation*}
differentiated yield:

\begin{alignat*}{9}
&dx &&=dt-ds \sin \vartheta \sin \zeta -s d\zeta \sin \vartheta \cos \zeta -s d \vartheta \sin \zeta \cos \vartheta, \\
&dy &&=\frac{dt}{\tan \zeta}-ds \sin \vartheta \cos \zeta +s d \zeta \sin \zeta \sin \vartheta - s d \vartheta \cos \zeta \cos \vartheta, \\
&dz &&=\frac{dt}{\sin \zeta \tan \vartheta}-ds \cos \vartheta +s d \vartheta \sin  \vartheta.
\end{alignat*}

\paragraph*{§27}

Before we proceed any further, it will not be out of place to have noted the principal relations among these formulas, and at first by eliminating $s$ we obtain these relations for the finite formulas:

\begin{alignat*}{9}
&x \cos \zeta - y \sin \zeta = t \cos \zeta -u \sin \zeta, \\
&x \sin \zeta +y \cos \zeta =t \sin \zeta +u \cos \zeta -s \sin \vartheta,\\
&x \sin \zeta \cos \vartheta +y \cos \zeta \cos \vartheta -z \sin \vartheta = t \sin \zeta \cos \vartheta +u \cos \zeta \cos \vartheta -v \sin \vartheta.
\end{alignat*}
Further, for the differentials we find the the following:

\begin{alignat*}{9}
& &&\text{I}^{\circ}. \quad && dx \cos \zeta -dy \sin \zeta =-s d \zeta \sin \vartheta, \\
& &&\text{II}^{\circ}. \quad && dx \sin \zeta +dy \cos \zeta =\frac{dt}{\sin \zeta}-ds \sin \vartheta - s d \vartheta \cos \vartheta \\
&\text{and} \quad && \text{III}^{\circ}. \quad && dx  \sin \zeta \cos \vartheta +dy \cos \zeta \cos \vartheta -dz \sin \vartheta =-s d \vartheta.
\end{alignat*}

\paragraph*{§28}

But because in this new calculation we reduced everything to the two variables $t$ and $s$, while in the first calculation we used the two variables $T$ and $U$, let us see, how these are expressed by those, and from the formulas found for $T$ and $U$ we indeed have

\begin{alignat*}{9}
& &&dT &&=\frac{dt \sin \omega}{\sin \zeta \sin \vartheta}&&- ds \sin \omega &&- s d\omega \cos \omega \\
&\text{and} \quad && dU &&=\frac{dt \cos \omega}{\sin \zeta \sin \vartheta}&&-ds \cos \omega &&+s d\omega \sin \omega;
\end{alignat*}
if we substitute these values in the formulas $dx$, $dy$ and $dz$ found before and carefully distinguish the two variables $t$ and $s$, we will obtain the following expression:

\begin{alignat*}{9}
&dx &&=dt \frac{l \sin \omega +\lambda \cos \omega}{\sin \zeta \sin \vartheta}&&-s d \omega (l \cos \omega -\lambda \sin \omega)&&-ds (l \sin \omega +\lambda \cos \omega),\\
&dy &&=dt \frac{m \sin \omega +\mu \cos \omega}{\sin \zeta \sin \vartheta}&&- s d\omega (m \cos \omega -\mu \sin \omega) &&-ds (m \sin \omega +\mu \cos \omega),\\
&dz &&=dt \frac{n \sin \omega +\nu \cos \omega}{\sin \zeta \sin \vartheta}&&- s d\omega (n \cos \omega -\nu \sin \omega) &&-ds (n \sin \omega +\nu \cos \omega),
\end{alignat*}
which we want to compare to those which arose in the last solution which are

\begin{alignat*}{9}
&dx &&=dt &&- s d \zeta \sin \vartheta \cos \zeta &&- s d \vartheta \sin \zeta \cos \vartheta &&- ds \sin \zeta \sin \vartheta,\\
&dy &&=\frac{dt}{\tan \zeta} &&+s d \zeta \sin \zeta \sin \vartheta &&- s d \vartheta \cos \zeta \cos \vartheta &&- ds \cos \zeta \sin \vartheta, \\
&dz &&=\frac{dt}{\sin \zeta \tan \vartheta} && &&+sd \vartheta \sin \vartheta &&-ds \cos \vartheta;
\end{alignat*}
and first, the terms affected by $ds$ have to be equal on both sides, whence we obtain these equations:

\begin{alignat*}{9}
&\text{I}^{\circ}. \quad && l\sin \omega &&+\lambda \cos \omega &&= \sin \zeta \sin \vartheta,\\
&\text{II}^{\circ}. \quad && m \sin \omega &&+\mu \cos \omega &&=\cos \zeta \sin \vartheta,\\
&\text{III}^{\circ}. \quad && n \sin \omega &&+\nu \cos \omega &&= \cos \vartheta.
\end{alignat*}

\paragraph*{§29}

Therefore, if  these values are now substituted in the first terms, which involve the differential $dt$ and those depending on it, namely $d \zeta$, $d \vartheta$ and $d\omega$, we will obtain the following equations:

\begin{alignat*}{9}
&l \cos \omega &&- \lambda \sin \omega &&= \frac{d \zeta \cos \zeta \sin \vartheta +d \vartheta \sin \zeta \cos \vartheta}{d \omega}&&=\frac{d(\sin \zeta \sin \vartheta)}{d \omega},\\
&m \cos \omega &&- \mu \sin \omega &&=\frac{-d \zeta \sin \zeta \sin \vartheta +d \vartheta \cos \zeta \cos \vartheta}{d \omega}&&=\frac{d(\cos \zeta \sin \vartheta)}{d \omega}, \\
&n \cos \omega &&-\nu \sin \omega &&=- \frac{d \vartheta \sin \vartheta}{d \omega}=\frac{d \cos \vartheta}{d \omega}.
\end{alignat*}
Here it is especially remarkable that  the one variable $s$ went out of these formulas found here completely such that now the quantities $l$, $\lambda$, $m$, $\mu$, $n$, $\nu$ are determined by the single variable $t$ and do not involve the other $s$ at all, whereas the quantities $T$ and $U$ contain both variables $t$ and $s$.

\paragraph*{§30}

Now, we found these two equations for defining the functions $l$ and $\lambda$:

\begin{alignat*}{9}
&l \cos \omega &&+\lambda \cos \omega &&= \sin \zeta \sin \vartheta,\\
&l \cos \omega &&-\lambda \sin \omega &&= \frac{d(\sin \zeta \sin \vartheta)}{d \omega}.
\end{alignat*}
Hence, the first equation multiplied by $\sin \omega +$ the second equation multiplied by $\cos \omega$ gives:

\begin{equation*}
l= \sin \zeta \sin  \vartheta \sin \omega +\cos \omega \frac{d(\sin \zeta \sin \vartheta)}{d\omega},
\end{equation*}
but I. $\cos \omega-$II. $\sin \omega$ gives:

\begin{equation*}
\lambda= \sin \zeta \sin \vartheta \cos \omega - \sin \omega \frac{d(\sin \zeta \sin \vartheta)}{d \omega}.
\end{equation*}
In the same way the remaining letters will be found as follows:

\begin{alignat*}{9}
&m &&= \cos \zeta \sin \vartheta \sin \medskip+\cos \omega \frac{d(\cos \zeta \sin \vartheta)}{d \omega},\\
& \mu &&= \cos \zeta \sin \vartheta \cos \omega - \sin \omega \frac{d(\cos \zeta \sin \vartheta)}{d \omega},\\
&n&&= \cos \vartheta \sin \omega +\frac{\cos \omega~ d \cos \vartheta}{d \omega},\\
&\nu &&= \cos \vartheta \cos \omega - \frac{\sin \omega~ d \cos \vartheta}{d \omega}.
\end{alignat*}
Behold these suitable values for the letters $l$, $\lambda$, $m$, $\mu$ and $n,\nu$ which are of such a nature that those three formulas $ldT+\lambda dU$, $mdT+\mu dT$ and $ndT+\nu dU$ become integrable and even the integrals  can easily be exhibited, which of course are

\begin{equation*}
x=t-s \sin \vartheta \sin \zeta, \quad y=u-s \sin \vartheta \cos \zeta, \quad z= v-s \cos \vartheta.
\end{equation*}
\paragraph*{§31}

Since our two solutions have to be completely identical to each other, there is no doubt that the remaining conditions mentioned above are also fulfilled, it will certainly be:

\begin{equation*}
ll+mm+nn=1, \quad \lambda \lambda +\mu \mu +\nu \nu=1, \quad l \lambda +m \mu +n \nu=0.
\end{equation*}
To show this, for the sake of brevity let us put

\begin{equation*}
\sin \zeta \sin \vartheta =p, \quad \cos \zeta \sin \vartheta=q \quad \text{and} \quad \cos \vartheta=r,
\end{equation*}
such that it is

\begin{equation*}
pp+qq+rr=1 \quad \text{and hence} \quad pdp+qdq+rdr=0,
\end{equation*}
now, because we have

\begin{alignat*}{9}
&l &&= p \sin \omega && +\frac{dp}{d \omega} \cos \omega, \\
&m &&= q \sin \omega &&+ \frac{dq}{d \omega} \cos \omega, \\
&n &&= r \sin \omega &&+ \frac{dq}{d \omega} \cos \omega, \\
&\lambda &&= p \cos \omega && -\frac{dp}{d \omega} \sin \omega, \\
& \mu &&=q \cos \omega &&- \frac{dq}{d\omega}\sin \omega, \\
&\nu &&=r \cos \omega && -\frac{dr}{d \omega}\sin \omega,
\end{alignat*}
having done the calculation we will hence find:

\begin{equation*}
1^{\circ}\quad ll+mm+nn=(pp+qq+rr)\sin^2 \omega +\frac{2 \sin \omega \cos \omega}{d \omega}(pdp+qdq+rdr)
\end{equation*}
\begin{equation*}
+\frac{\cos^2 \omega}{d \omega^2}(dp^2+dq^2+dr^2)
\end{equation*}
or

\begin{equation*}
ll+mm+nn=\sin^2 \omega + \frac{\cos^2 \omega}{d \omega^2}(dp^2+dq^2+dr^2),
\end{equation*}
and so the whole question is now shifted to the investigation of the value $dp^2+dq^2+dr^2$. But because it is

\begin{alignat*}{9}
& &&dp &&= + d \zeta \cos \zeta \sin \vartheta &&+d \vartheta \sin  \zeta \cos \vartheta, \\
& &&dq &&=- d \zeta \sin \zeta \sin \vartheta &&+ d \vartheta \cos \zeta \cos \vartheta \\
&\text{and} \quad && dr&&=-d \vartheta \sin \vartheta,
\end{alignat*}
we conclude

\begin{equation*}
dp^2+dq^2+dr^2= d \zeta^2 \sin^2 \vartheta +d \vartheta^2 =d \omega^2,
\end{equation*}
so that it is certain that it is

\begin{equation*}
\frac{dp^2+dq^2+dr^2}{d \omega^2}=1,
\end{equation*}
whence it is manifest that it will be:

\begin{equation*}
ll+mm+nn= \sin^2 \omega +\cos^2 \omega =1.
\end{equation*}

\paragraph*{§32}

In like manner  we will find for the Greek letters:

\begin{equation*}
\lambda \lambda +\mu \mu +\nu \nu = (pp+qq+rr) \cos^2 \omega -\frac{2 \sin \omega \cos \omega}{d \omega}(pdp+qdq+rdr)
\end{equation*}
\begin{equation*}
+\frac{\sin^2 \omega}{d \omega^2}(dp^2+dq^2+dr^2),
\end{equation*}
which manifestly, as before, yields 

\begin{equation*}
\lambda \lambda +\mu \mu +\nu \nu =\cos^2 \omega +\sin^2 \omega =1.
\end{equation*}
Therefore, it remains to examine the third property, for which we obtain:

\begin{alignat*}{9}
&l \lambda &&pp \sin \omega \cos \omega &&-\frac{pdp}{d \omega}\sin^2 \omega &&+\frac{pdp}{d \omega}\cos^2 \omega &&- \frac{dp^2}{d \omega^2} \sin \omega \cos \omega, \\
&m \mu &&qq  \sin \omega \cos \omega &&-\frac{qdq}{d \omega}\sin^2 \omega &&+\frac{qdq}{d \omega}\cos^2 \omega &&- \frac{dq^2}{d \omega^2} \sin \omega \cos \omega, \\
&n \nu &&rr  \sin \omega \cos \omega &&-\frac{rdr}{d \omega}\sin^2 \omega &&+\frac{rdr}{d \omega}\cos^2 \omega &&- \frac{dr^2}{d \omega^2} \sin \omega \cos \omega;
\end{alignat*}
having collected these into  one sum, it will be

\begin{equation*}
l \lambda +m \mu +n \nu = \sin \omega \cos \omega -\sin \omega \cos \omega=0.
\end{equation*}
And  this way we solved the analytical problem mentioned above (§ 7), which solution can be given in short form as follows. \newpage

\subsection*{Analytical Problem}

\paragraph*{§33}
\textit{Given two variables $T$ and $U$, to find six functions $l$, $m$, $n$ and $\lambda$, $\mu$, $\nu$  of such a nature that the following six conditions are satisfied:}

\begin{alignat*}{9}
&\text{I}^{\circ}. ~~ \left(\frac{dl}{dU}\right)=\left(\frac{d \lambda}{dT}\right), \quad \text{II}^{\circ}.~~ \left(\frac{dm}{dU}\right)=\left(\frac{d\mu}{dT}\right),\quad \text{III}^{\circ}.~~ \left(\frac{dn}{dU}\right)=\left(\frac{d \nu}{dT}\right),\\
&\text{IV}^{\circ}~~ ll+mm+nn=1,\quad \text{V}^{\circ}.~~ \lambda \lambda+\mu \mu +\nu \nu =1,\\
&\text{VI}^{\circ}.~~ l \lambda +m \mu +n \nu=0.
\end{alignat*}

\subsection*{Solution}

Having introduced the two new variables $s$ and $t$ into the calculation, introduce  two functions $\zeta$ and $\vartheta$ of the letter $t$, which functions are to be considered as angles, of course;  from these form the angle  $\omega$ in such a way that it is

\begin{equation*}
d \omega =\sqrt{d \zeta^2 \sin^2 \vartheta +d \vartheta^2}.
\end{equation*}
Then from this the two variables $T$ and $U$ are indeed determined in such a way that it is

\begin{alignat*}{9}
&T &&= \int \frac{dt \sin \omega}{\sin \zeta \sin \vartheta}&&- s \sin \omega, \\
&U&&= \int \frac{dt \cos \omega}{\sin \zeta \sin \vartheta}&&-s \cos \omega;
\end{alignat*}
having done this the six functions in question will behave as this

\begin{alignat*}{9}
&l &&= \sin \zeta \sin \vartheta \sin \omega +\frac{\cos \omega}{d\omega}d(\sin \zeta \sin \vartheta),\\
&\lambda &&= \sin \zeta \sin \vartheta \cos \omega -\frac{\sin \omega}{d\omega}d(\sin \zeta \sin \vartheta),\\
&m &&= \cos \zeta \sin \vartheta \sin \omega +\frac{\cos \omega}{d\omega}d(\cos \zeta \sin \vartheta),\\
&\mu &&= \cos \zeta \sin \vartheta \cos \omega -\frac{\sin \omega}{d\omega}d(\cos \zeta \sin \vartheta),\\
&n &&= \cos \vartheta \sin \omega +\frac{\cos \omega}{d \omega}d \cos \vartheta, \\
&\nu &&= \cos \vartheta \cos \omega -\frac{\sin \omega }{d \omega}d \cos \vartheta.
\end{alignat*}
But by means of these three values the following three differential formulas:

\begin{equation*}
\text{I}^{\circ}. ~~ ldT+\lambda dU, \quad \text{II}^{\circ}.~~ mdT+\mu dU, \quad \text{III}^{\circ}.~~ ndT+\nu dU,
\end{equation*}
in which the first three conditions are contained, of course, are not only rendered integrable, but also the integrals will be expressed as follows:

\begin{alignat*}{9}
&\text{I}^{\circ}. \quad &&\int(ldT &&+\lambda dU)&&=t-s \sin \vartheta \sin \zeta,\\
&\text{II.}^{\circ}. \quad && \int (mdT&&+ \mu dU)&&= \int \frac{dt}{\tan \zeta}- s \sin \vartheta \cos \zeta, \\
&\text{III.}^{\circ}. \quad && \int(ndT&&+\nu dU)&&=\int \frac{dt}{\sin \zeta \tan \vartheta}-s \cos \vartheta,
\end{alignat*}
which solution is therefore to be considered as complete because it contains two arbitrary functions.

\paragraph*{§34}

This expansion without any doubt is of greatest importance and especially deserves it that we  inquire its single elements with all eagerness. And at first, because having introduced the letters $p$, $q$ and $r$ in such a way that it is

\begin{equation*}
pp+qq+rr=1, \quad \text{and} \quad dp^2+dq^2+dr^2=d \omega^2,
\end{equation*}
we found

\begin{equation*}
l \sin \omega +\lambda \cos \omega =p \quad \text{and} \quad l \cos \omega - \lambda \sin \omega =\frac{dp}{d \omega},
\end{equation*}
if we differentiate the first equation, we will have

\begin{equation*}
d l \sin \omega +d \lambda \cos \omega + l d\omega \cos \omega - \lambda d \omega \sin \omega =dp
\end{equation*}
and hence

\begin{equation*}
dl \sin \omega +d \lambda \cos \omega=0,
\end{equation*}
such that it is

\begin{equation*}
\frac{d \lambda}{d l}=- \tan \omega.
\end{equation*}
In like manner we will indeed also find

\begin{equation*}
\frac{d\mu}{dm}=-\tan \omega \quad \text{and} \quad \frac{d \nu}{dn}=- \tan \omega.
\end{equation*}
Therefore, behold this most beautiful property which intercedes among our six functions $l$, $m$, $n$ and $\lambda$, $\mu$, $\nu$ and which can also be represented in this way that it is

\begin{equation*}
dl:d \lambda =dm:d \mu=dn:d \nu =- \cos \omega: \sin \omega.
\end{equation*}

\paragraph*{§35}

Therefore, if we  consider those things carefully, we will discover certain traces and following them we will be able to find a direct solution of this most difficult problem. Of course, having constituted these equations:

\begin{equation*}
dx=ldT+\lambda dU, \quad dy=mdT+\mu dU, \quad dz=ndT+\nu dU
\end{equation*}
it is convenient to observe at first that the quantities $l$, $m$, $n$ and $\lambda$, $\mu$, $\nu$ have to be functions of one single new variable, which nevertheless has a certain relation to the two principal variables $T$ and $U$. Therefore, let $\omega$ be this new variable and let also our six quantities  be certain functions of it. And we have already seen, if the letters $p$, $q$ and $r$ are such functions of $\omega$, that it  then  is

\begin{equation*}
pp+qq+rr=1 \quad \text{and} \quad dp^2+dq^2+dr^2=d\omega^2,
\end{equation*}
that then by putting:

\begin{alignat*}{9}
& l &&= p \sin \omega &&+\frac{dp}{d \omega} \cos \omega, \\
&m  &&=q \sin \omega &&+ \frac{dq}{d \omega}\cos \omega, \\
&n &&= r \sin \omega && +\frac{dr}{d \omega} \cos \omega, \\
&\lambda &&= p \cos \omega &&- \frac{dp}{d \omega}\sin \omega, \\
& \mu &&=q \cos \omega && - \frac{dq}{d\omega}\sin \omega, \\
&\nu &&= r \cos \omega && -\frac{dr}{d \omega}\sin \omega,
\end{alignat*}
now these three conditions are already fulfilled, namely:

\begin{equation*}
ll+mm+nn=1,\quad \lambda \lambda +\mu \mu +\nu \nu=1 \quad \text{and} \quad l \lambda +m \mu +n \nu=0;
\end{equation*}
furthermore, from this we  already deduced  the extraordinary property that it is

\begin{equation*}
d \lambda=-dl \tan \omega, \quad d \mu =-dm \tan \omega \quad \text{and} \quad -dn \tan \omega,
\end{equation*} 
which will be of immense use for us to fulfill the remaining conditions, as it will become clear soon.

\paragraph*{§36}

These three conditions certainly demand that those differential formulas exhibited for $dx$, $dy$ and $dz$ are rendered integrable; for this one has to find the relation among the two variables $T$ and $U$ and  $\omega$. To achieve this, by integrating convert these differential equations into the following forms:

\begin{alignat*}{9}
&x &&=lT &&+\lambda U &&- \int (Tdl &&+U d\lambda),\\
&y &&=mT &&+\mu U &&- \int (Tdm &&+U d \mu), \\
&z &&=nT &&+\nu U &&- \int (Tdn &&+U d \nu);
\end{alignat*}
now, these three new integral formulas indeed will obtain the following forms:

\begin{alignat*}{9}
&x&&=l T &&+\lambda U &&- \int dl &&(T-U \tan \omega), \\
&y&&=m T &&+\mu U &&- \int dm &&(T-U \tan \omega), \\
&z&&=n T &&+\nu U &&- \int dn &&(T-U \tan \omega).
\end{alignat*}
Since $l$,$m$, $n$ are functions of the same variable $\omega$, it is manifest that these three formulas are indeed rendered integrable, if  the expression $T-U \tan \omega$ was a function of the new variable $\omega$; hence, if such a function is indicated by the letter $\Omega$, we will have

\begin{equation*}
T- U \tan \omega =\Omega;
\end{equation*}
by means of these equations the equation in question between the variables $T$, $U$ and $\omega$ is determined.

\paragraph*{§37}

Hence, if for $\Omega$  an arbitrary function of $\omega$ is taken, of which also, as we saw, the letters $p$, $q$ and $r$ are certain functions, by means of which we already defined the letters $l$, $m$, $n$ and $\lambda$, $\mu$, $\nu$, the two variables $T$ and $U$ must be of such a nature that it is $T=\Omega +U \tan \omega$; of course, we  want to keep only the two variables $U$ and $\omega$ in the calculation and therefore let us introduce this value instead of $T$, then our three integral formulas can be represented this way:

\begin{alignat*}{9}
&x&&= l \Omega &&+lU \tan \omega && + \lambda U &&- \int \Omega dl, \\
&y&&=m \Omega &&+m U \tan \omega &&+ \mu U &&- \int \Omega dm, \\
&z&&=n \Omega &&+n U \tan \omega && +\nu U &&- \int \Omega dn, 
\end{alignat*}
which expressions are easily transformed into the following ones

\begin{alignat*}{9}
&x &&= U(l \tan \omega &&+\lambda)&&+\int ld \Omega &&= \frac{Up}{\cos \omega}&&+\int p \sin \omega d \Omega &&+\int \frac{dp d \Omega}{d \omega}\cos \omega,\\
&y &&= U(m \tan \omega &&+\mu)&&+\int md \Omega &&= \frac{Uq}{\cos \omega}&&+\int q \sin \omega d \Omega &&+\int \frac{dq d \Omega}{d \omega}\cos \omega,\\
&z &&= U(n \tan \omega &&+\nu)&&+\int nd \Omega &&= \frac{Ur}{\cos \omega}&&+\int r \sin \omega d \Omega &&+\int \frac{dr d \Omega}{d \omega}\cos \omega.
\end{alignat*}

\newpage

\subsection*{Third solution of the principal problem derived from the theory of light and shadow}

\paragraph*{§38}

What is usually treated in Optics on light and shadow, is mostly restricted to the highly special case in which both the shining and the opaque body, from which the shadow is projected, have a spherical shape; hence either a cylindrical or a conic or a convergent or divergent shadow arises, depending on whether the opaque body was either equal to or  smaller or larger  than the shining body. But whenever the shape  of either the shining or the opaque body or of both recedes from that of a sphere, we hardly find anything, which could  content us, in the books written on this subject; if we wanted to treat this subject in  general, attributing any shapes to both bodies, the shining and the opaque, a most difficult question would arise; and this question belongs to that part of the analysis of the Infinite on functions of two or more variables which was begun  to be constructed not so long ago.

\paragraph*{§39}

But what especially  extends to our undertaking from this theory, is that the shapes of the shadows are always of such a nature that their surface can be unfolded onto a plane; hence it is vice versa understood, if we were able to determine the shape of the shadow for any figure of both the shining and the opaque body, that then at the same time also our problem will be perfectly solved.

\paragraph*{§40}

That the shape of the shadow is indeed always subjected to our problem can easily be shown this way. Since the shadow is terminated by the most outer rays of the shining body which at the same time touch upon the opaque body, first it is plain that  infinitely many straight lines are given in the surface of a certain shadow, since the single rays proceed in straight lines; furthermore, all these rays will touch upon both the shining and the opaque body, whence, if  any plane is imagined, which those two bodies touch at the same time and the point of contact on the shining body is denoted by the letter $M$, on the opaque on the other hand by the letter $m$, it is clear that the straight line $Mm$, if it is elongated, exhibits the ray of the light by which the shadow is terminated, which is also to be understood about the other infinitely close rays, which are emitted from the point $M$ on the same tangential plane, which rays can also be considered as the tangents of the opaque body, from which the most excellent property of our problem arises that any two infinitely close lines to be drawn in the surface at the same time are found in the same plane.

\paragraph*{§41}

But this theory of light and shadow extends too far to be discussed here in more detail; therefore, we will only take everything necessary to solve our  present problem from it. Having put  the shape of both the shining and the opaque body aside let us investigate only the shape of a shadowy cone; for this purpose, let us consider two parallel sections removed from a each other by a given distance. Because it is possible to attribute any arbitrary shape to these sections, it is manifest that this consideration contains completely all shapes of the shadows.

\paragraph*{§42}

Therefore, let (Fig. 6) these two sections be normal to the plane of the diagram and be based perpendicularly on the line $Aa$,
\begin{center}
\includegraphics[scale=0.4]{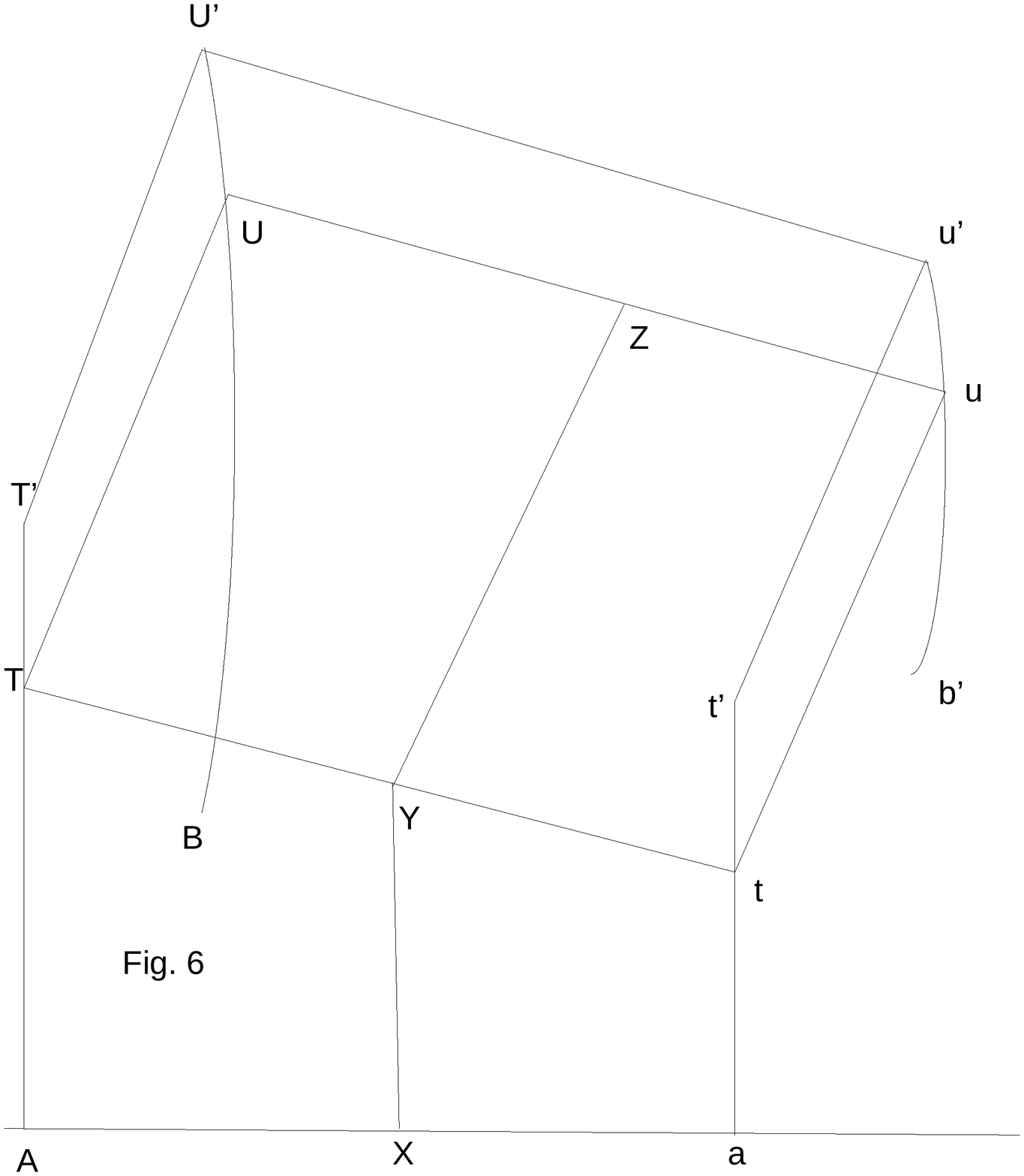}
\end{center}
 and  at first let $BUU^{\prime}$ be the curve, whose nature  is to be expressed by an equation between the coordinates $AT=T$ and $TU=U$; in the same way let $buu^{\prime}$ be another curve different from the first, for which an equation between the coordinates $at=t$, $tu=u$  is given, but put the interval between these sections $Aa=a$; here, it will certainly be possible to consider the one section $BUU^{\prime}$ as a shining planar disk, and while the other $buu^{\prime}$ refers to  a opaque planar disk, the shadow-cone, we contemplate,  will  arise from the light rays.

\paragraph*{§43}

But let the points $U$ and $u$ be taken in such a way that the line $Uu$, if it is elongated, represents  the ray terminating the shadow; because it  has to lie in the plane touching both  disks, it is necessary that both elements $UU^{\prime}$ and $uu^{\prime}$ lie in the same plane together with the line $Uu$. From this it is clear that these two elements are parallel to each other, whence it follows that  the ratio between the differentials has to be the same such that it is $dT:dU=dt:du$; this is why, if it one puts $dU=\varphi dT$, it will also be $du= \varphi dt.$

\paragraph*{§44}

Therefore, consider this quantity $\varphi$ as the principal variable, by which all the remaining ones are determined in the following way. For the first curve $BU$ let $T$ be a function of $\varphi$, whose nature defines the properties of the curve $BUU^{\prime}$, but then it will be

\begin{equation*}
dU= \varphi dT \quad \text{and} \quad U= \int \varphi dT;
\end{equation*}
it is evident that  this way any arbitrary curve can be expressed by means  of the variable $\varphi$. In like manner for the other curve $buu^{\prime}$ the abscissa $t$ will certainly become equal to a function of $\varphi$ and then one will equally have

\begin{equation*}
du=\varphi dt \quad \text{and} \quad u= \int \varphi dt, 
\end{equation*}
whence, because the two curves are completely arbitrary, it is possible to assume any functions of $\varphi$ for the letters $T$ and $t$; having done so at the same time the two ordinates $U$ and $u$ are determined.

\paragraph*{§45}

Now, let us take an arbitrary point $Z$ on the line $Uu$; because this point lies on the surface we investigate, let us drop the perpendicular $ZY$ intersecting the line $Tt$ to the plane of the diagram from that point and from $Y$ let us draw the normal $YX$ to our axis $Aa$ that for the indefinite point $Z$ we obtain three coordinates we want to call:

\begin{equation*}
AX=x, \quad XY=y \quad \text{and} \quad YZ=z,
\end{equation*}
and now it will be easy to find an equation between these three coordinates by means of which the nature of the surface in question is expressed.

\paragraph*{§46}

The principles of geometry immediately give us these properties

\begin{alignat*}{9}
&T-&&t:a&&=T-y:x, \quad && \text{or} \quad &&Tx &&-tx &&=aT &&-ay, \\
&U-&&u:a&&=T-z:x, \quad && \text{or} \quad &&Ux &&-ux &&=aU &&-az, 
\end{alignat*}
whence by means of the two variables $\varphi$ and $x$ it will be possible to define the two coordinates $y$ and $z$, since we will have:

\begin{equation*}
y=T-\frac{x(T-t)}{a}\quad \text{and} \quad z=U-\frac{x(U-u)}{a};
\end{equation*}
for, if  the variable $\varphi$ together with the ones depending on it $T$, $t$ and $U$, $u$ is eliminated from these two equations, an equation expressing the nature of our surface will result.

\paragraph*{§47}

But we do not want to do such an elimination explicitly, since the nature of the surface can be seen a lot more  clearly from the two  equations we found, which  are already so simple that it would be a crime to desire a more convenient solution; but it will nevertheless be useful to manipulate the forms of these equations a little bit. In a more general way let us represent the values for $y$ and $z$ as this

\begin{equation*}
y=P+Qx \quad \text{and} \quad z=R+Sx,
\end{equation*}
where the letters $P$, $Q$, $R$, $S$ now denote functions of the other variable $\varphi$, and now the question is: Of what nature must these functions be that the two exhibited equations define a surface  which can be unfolded onto a plane?

\paragraph*{§48}

Therefore, let us compare these assumed forms to those we found before and we will have:

\begin{equation*}
P=T \quad \text{and} \quad R=U, \quad Q=\frac{t-T}{a},\quad S=\frac{u-U}{a};
\end{equation*}
here, because  $T$ and $t$ are arbitrary functions of $\varphi$, it is evident that the functions $P$ and $Q$ can be taken arbitrarily, and since $U$ and $u$ depend on $T$ and $t$, the functions $R$ and $S$ will also have to depend on the first two $P$ and $Q$ in a certain way. But because it is

\begin{equation*}
T=P, \quad t=P+aQ, \quad U=R \quad \text{and} \quad u=R+aS,
\end{equation*} 
let us substitute these values in the fundamental formulas

\begin{equation*}
dU=\varphi dT \quad \text{and} \quad du =\varphi dt
\end{equation*}
and we will obtain

\begin{equation*}
dR=\varphi dP \quad \text{and} \quad dR+adS=\varphi dP+a \varphi dQ
\end{equation*}
or $dS=\varphi dQ$.

\paragraph*{§49}

Therefore, we will also be able to eliminate the quantity $\varphi$ from the calculation, because it either is $\varphi=\frac{dR}{dP}$ or $\varphi = \frac{dS}{dQ}$, such that instead of it one of the letters $R$ and $S$ are now arbitrary;  hence, if  $P$, $Q$ and $R$ were any arbitrary functions of a certain variable, which is the same for all three functions\footnote{and not necessarily $\varphi$.}, then $S$ has to be such a function of the same variable that:

\begin{equation*}
dS= \frac{dQdR}{dP} \quad \text{or} \quad \frac{dS}{dR}=\frac{dQ}{dP};
\end{equation*}
this solution can even be rendered more beautiful in such a way that we say that  for the letters $P$, $Q$, $R$, $S$ one has to assume functions of certain variable of such a kind that it is $\frac{dS}{dR}=\frac{dQ}{dP}$ or even $\frac{dS}{dQ}= \frac{dR}{dP}$; having done so these two equations

\begin{equation*}
y=P+Qx \quad \text{and} \quad z=R+Sx
\end{equation*}
will express the nature of the solid in question.

\paragraph*{§50}

It does not matter by which letter the variable, of which $P$, $Q$, $R$ and $S$ are functions, is indicated; one can even take one of these $P$, $Q$, $R$, $S$ for it. The remaining ones are then to be understood as functions of that variable. Therefore, as long as one of them retains a constant value, the remaining ones will also be constant, and then from the variability of $x$ all straight lines  which can be drawn on the surface will arise.

\paragraph*{§51}

The prescribed condition $\frac{dS}{dQ}=\frac{dR}{dP}$ will manifestly be fulfilled by taking the quantities $P$ and $R$ as constants; hence, a particular solution of our problems follows. For, let us assume that it is $P=A$ and $R=B$, such that now $S$ is to be considered as a function of $Q$. But it is always possible to vary the coordinates in such a way that $A=0$ and $B=0$; having done this, because of $Q=\frac{y}{x}$, $\frac{z}{x}=S$ will be  a homogeneous function of no dimension of $x$ and $y$, or $z$ will become equal to a homogeneous function of one dimension of $x$ and $y$ which is the criterion for a surface of a cone. 

\paragraph*{§52}

The condition is also fulfilled by taking $Q=0$ and $S=0$ such that $R$ remains a function of $P$, in which case  a function of $y$ will arise for $z$; because the final equation involves only two variables, $y$ and $z$, it will be an equation for a cylindrical solid;  the same happens, if we put either $P=0$ and $Q=0$ or $R=0$ and $S=0$; for, in the first case one has $y=0$, in the second on the other hand $z=0$, in both cases it is the equation for a plane.

\paragraph*{§53}

But to understand  also other species of solids of this kind and to find the simpler cases let us assume :

\begin{equation*}
P=a\varphi^{\alpha}, \quad Q=b\varphi^{\beta},\quad R=c \varphi^{\gamma}, \quad S=d \varphi^{\delta},
\end{equation*}
and to fulfill the prescribed condition it is necessary that it is

\begin{equation*}
\frac{b \beta}{a \alpha}\varphi^{\beta -\alpha}=\frac{d \delta}{c \gamma}\varphi^{\delta -\gamma},
\end{equation*}
whence two conditions arise, the first for the exponents

\begin{equation*}
\beta -\alpha = \delta -\gamma,
\end{equation*}
the other for the coefficients:

\begin{equation*}
\frac{b \beta}{a \alpha}= \frac{d \delta}{c \gamma},
\end{equation*}
both of which conditions are met by taking the values as follows:

\begin{equation*}
a=\frac{fg}{\varkappa +\lambda}, \quad b=\frac{fh}{\varkappa +\mu}, \quad c=\frac{gk}{\lambda +\nu}, \quad d=\frac{hk}{\mu+\nu},
\end{equation*}
\begin{equation*}
\alpha = \varkappa + \lambda, \quad \beta = \varkappa +\mu, \quad \gamma= \lambda +\nu, \quad \delta= \mu +\nu,
\end{equation*}
then the equations will be:

\begin{equation*}
y= a \varphi^{\alpha}+b\varphi^{\beta}x, \quad z=c \varphi^{\gamma}+d \varphi^{\delta} x.
\end{equation*}

\paragraph*{§54}

Therefore, substituting some explicit numbers for the powers  let us consider this example:

\begin{equation*}
y= 2 \varphi +3 \varphi^2 x \quad \text{and} \quad z=\varphi^2 +2 \varphi^3 x,
\end{equation*}
whence after the elimination of the letter $\varphi$ the following equation is found\footnote{In Euler's original paper he gives the equation: $4y^3x+72y^2xxz-yy-18yxz+27xxzz+2z=0$; this equation does not give a developable surface. The following equation is correct as  pointed out by A. Speiser in the Opera Omnia version.}

\begin{equation*}
-4xy^3-y^2+18xyz+27x^2z^2+4z=0,
\end{equation*}
which is therefore for a solid whose surface can be unfolded onto the plane.

\end{document}